\documentclass[12pt]{article}
\usepackage{amsmath}
\usepackage{enumitem}
\usepackage{authblk}
\usepackage{lineno}
\usepackage[toc,page]{appendix}
\usepackage{amsthm}
\usepackage{changepage}
\usepackage{aliascnt}
\usepackage{theoremref,amsthm}
\usepackage{morefloats}
\usepackage{amsfonts}
\usepackage{amssymb}
\usepackage{graphics}
\usepackage{graphicx}
\usepackage{listings}
\usepackage{xcolor}
\usepackage{float}
\usepackage{verbatim}
\usepackage{geometry}
\usepackage{centernot}
\usepackage[utf8]{inputenc}
\usepackage[english]{babel}
\usepackage{setspace}
\usepackage{thmtools}
\usepackage{xcolor}
\usepackage[bottom]{footmisc}
\usepackage{euscript}
\usepackage{mathtools}
\DeclareMathAlphabet{\itbf}{OML}{cmm}{b}{it}

\newcommand{\bea}{\begin{eqnarray*}}
\newcommand{\eea}{\end{eqnarray*}}
\newcommand{\bean}{\begin{eqnarray}}
\newcommand{\eean}{\end{eqnarray}}
\newcommand{\p}{\partial}
\newcommand{\f}{\frac}
\newcommand{\s}{\sqrt}
\newcommand{\ds}{\displaystyle}
\newcommand{\no}{\nonumber}
\newcommand\ov{\overline}
\newcommand{\ri}{\rightarrow}
\newcommand{\sm}{\setminus}
\newcommand{\aaa}{\mbox{$[$}}
\newcommand{\bbb}{\mbox{$]$}}

\newtheorem{lem}{Lemma}[section]

\newtheorem{thm}{Theorem}[section]

\newcommand{\doubleR}{\mathbb{R}}
\newcommand{\RR}{\mathbb{R}}

\renewcommand{\a}[1]{\left\vert #1 \right\vert}
\newcommand{\norm}[1]{\left\lVert#1\right\rVert}

\renewcommand{\ov}[1]{\overline{#1}}
\newcommand{\pfrac}[2]{\frac{\partial #1}{\partial #2}}

\newcommand{\bm}{{\itbf m}}
\newcommand{\bp}{{\itbf p}}
\newcommand{\bt}{{\itbf t}}
\newcommand{\bq}{{\itbf q}}

\newcommand{\bx}{{\itbf x}}
\newcommand{\bv}{{\itbf v}}

\newcommand{\bh}{{\itbf h}}
\newcommand{\bev}{{\itbf e}}

\newcommand{\bn}{{\itbf n}}
\newcommand{\by}{{\itbf y}}

 \author{ Darko  Volkov and  Yulong Jiang
\thanks{Department of Mathematical Sciences,
Worcester Polytechnic Institute, Worcester, MA 01609. Corresponding author email: darko@wpi.edu. 
} }

\begin{document}
		\title{Stability properties of a crack inverse problem in half space}
		\maketitle

\begin{abstract}
We show in this paper  a Lipschitz stability 
result for a crack
inverse problem in half space. The direct problem is a Laplace equation
with zero Neumann condition on the top boundary.
The forcing term is a discontinuity across the crack.
This formulation can be related to geological faults in elastic media or to 
irrotational incompressible flows in a half space minus an inner wall.
The  direct problem is well posed in an appropriate 
functional space. We study the related inverse problem where
the jump across the crack is unknown, and more importantly, the geometry and the 
location of the crack are unknown.  
The data for the inverse problem is of Dirichlet type over a portion of the top boundary.
We prove that this inverse problem is uniquely solvable 
under some assumptions on the geometry of the crack.
The highlight of this paper is showing a stability result for this inverse problem.
Assuming that the crack is planar, we show that reconstructing 
the plane containing the crack is Lipschitz stable despite the fact that the forcing term for the underlying PDE is unknown. 
This uniform stability result holds under the assumption that the forcing
term is bounded above and the Dirichlet data is bounded below away from zero
in appropriate norms.  
\end{abstract}

\section{Introduction}
We consider a 
 PDE model  derived from a problem in geophysics 
where seismic or displacement data is collected by surface sensors and then processed
to reconstruct a source and  an unknown fault.
Recently, there has been some progress
in the 
 understanding of the 
underlying mathematics.
In particular, well posedness of the forward problem for 
the three dimensional linear elasticity model
and uniqueness
for the related inverse problem were shown in    \cite{volkov2017determining}, 
and in \cite{aspri2020analysis, aspri2020arxiv} where more general fault geometries
and elasticity tensors were studied.
A stability result in the case of planar faults was achieved in
  \cite{triki2019stability}.\\
	Here, we examine  the case of a model involving the Laplace equation.
	This model is also relevant to geophysics: in  dimension two,
	it relates to the so called anti-plane strain configuration.
	This  configuration has already attracted much attention 
	from geophysicists and mathematicians due to 
	how simple and yet relevant
	  the formulation is
	\cite{ionescu2008earth, ionescu2006inverse, dascalu2000fault}.
 In dimension three, the model considered in this paper can be related to irrotational incompressible 
flows in a medium with a top wall and an inner wall.\\
	We now outline the contents of this paper.
	In section 
	\ref{direct and inverse crack problems},
	we first write out the governing partial differential equation for this problem 
	for which we  prove existence and uniqueness in an adequate functional space.
	This equation involves the Laplace operator for the unknown $u$ in the lower half space
	minus a crack $\Gamma$.  Mathematically,  $\Gamma$ is set to be 
	an open Lipschitz surface.
	A zero Neumann condition is prescribed on the top boundary ${x_3 =0}$ while 
	the jump of $u$ across $\Gamma$ is set to be equal to some $g$.
	Next, we introduce the related crack inverse problem. 
	For that problem, 
	 the additional knowledge of Dirichlet type  data $u |_V$
	on a relatively open subset $V$ of the top boundary ${x_3 =0}$ is provided.
	The challenge is that the geometry and the location of $\Gamma$
	are unknown,
	as well as  the jump $g$.
	Assuming that the crack $\Gamma$ is a bounded 
	Lipschitz graph of a function in $ (x_1,x_2)$
	defined on a relatively open subset of the top boundary and one of the three conditions in
	Theorem \ref{InverseProblemResult} is satisfied, we are able to prove that
	the  crack inverse problem is uniquely solvable.  If all three conditions fail to be satisfied
	the  crack inverse problem may not be uniquely solvable: we give an example in Appendix
	\ref{counter-example}.
	\\
	Even in the simple case where $\Gamma$ is known,   mapping
	 $g $  to $u |_V$  is a linear and compact operation, so its inverse is  unbounded.
	To obtain a stability result, we assume that $\Gamma$ is included in the plane
	with equation $x_3 = a x_1 + b x_2 +d$. The stability result, which is the highlight of this paper,
	is formulated 
	in the geometry parameter $\bm= (a,b,d) $.
	Theorem  \ref{main theorem},
	states that  $\bm $
	depends on $u |_V$  in a Lipschitz fashion 
	for all jumps $g$  bounded above in $H^1_0$ and such that the corresponding 
	$L^2$ norm of $u |_V$  remains bounded below by a positive constant.\\
	To put our work in perspective, let us first explain how our results
 compare to 
	\cite{friedman1989determining}. 
	In that paper, it was shown in a two dimensional setting that
applying two
adequately chosen forcing terms on the boundary can uniquely determine cracks.
 The proof technique used in that paper
is entirely different from the technique used here which is based on the inverse function theorem 
and forming a PDE on the surface  $\Gamma$. 
The inverse problem considered 
in \cite{friedman1989determining} is active, which means  that boundary conditions 
can be  adequately chosen for the result to hold. 
In contrast,   we are  interested in a passive inverse problem, so the forcing term
$g$ is an unknown for the inverse problem. 
More recently, the authors of  \cite{aspri2020analysis} studied  
a fault inverse problem closely related to the crack inverse problem considered here.
That problem was formulated using 
the three dimensional elasticity equations 
 and was analyzed
from the viewpoint of very weak solutions. 
Unique-continuation arguments were used to establish uniqueness for  the 
crack inverse problem.
In \cite{aspri2020arxiv}
the case of layered media with 
piecewise Lipschitz Lam\'e coefficients was covered:
uniqueness for   the crack inverse problem 
was established under minimal assumptions for the fault geometry.
%
In \cite{triki2019stability},
a stability result closely related to Theorem \ref{main theorem} 
of this present paper was proved 
for the  linear elasticity equation  case.
Due to the higher level of complexity of the elasticity equations, additional 
assumptions had to be made in order to obtain 
a stability result.  The jump  $g$ across $\Gamma$ in that case was 
 modeled to be a 
 tangential vector field. To obtain a  stability estimate, it was assumed 
that $\Gamma$ is not horizontal,
that  $g$ is one-directional,
or that  $g$ is  the gradient of a function in $H^2$.
Theorem 4.1 in \cite{triki2019stability} provides a
 stability estimate which is uniform in all  planar faults $\Gamma$ 
and jumps $g$, provided a geometry $\Gamma'$ and a jump $g'$ are fixed.
In this paper, in Theorem \ref{main theorem}, our uniform Lispchitz statement is 
stronger because it is also valid for all $\Gamma'$ and all
$g'$  bounded above and below in appropriate norms.

\section{The direct and inverse crack problems in half space}\label{direct and inverse crack problems}
\subsection{Problem formulation and physical interpretation}
Let $\RR^{3-}$ be the open half space  $\{ x_3<0 \}$, where
we  use the notation $\bx=(x_1,x_2,x_3)$ for vectors in $\RR^3$.
  Let $\Gamma$ be a Lipschitz  open surface in $\RR^{3-}$ and $D$ a domain in $\RR^{3-}$
	with Lipschitz boundary such that $\Gamma \subset \p D$.
	We assume that $\Gamma$ is strictly included in $\RR^{3-}$ so that  the distance 
	from $\Gamma$ to the plane $\{ x_3 =0 \}$ is positive.
We define the direct crack problem to be the boundary value problem, 
	\bean
        \Delta u=0\text{ in }\doubleR^{3-}\sm \ov{\Gamma},  \label{BVP1}     \\
        \p_{x_3} u=0\text{ on  the surface } x_3=0,  \label{BVP2}     \\
     \aaa \frac{\p u}{\p \bn} \bbb = 0 \mbox{ across }\Gamma,   \label{BVP3}\\
      \aaa u \bbb
			=g \mbox{ across }\Gamma,  \label{BVP4}\\
    u(\bx)=O ({\frac{1}{| \bx|}}) \text{ uniformly as } |\bx|\to \infty,     \label{Decay1}
\eean
where  $\aaa v \bbb $ denotes the jump of a function $v$ across $\Gamma$ in the normal
direction,
and $\bn$ is a unit normal vector  to $\Gamma$.
In problem (\ref{BVP1}-\ref{Decay1}),
$u$ can be a model for the potential of an irrotational flow,
with an impermeable and immobile wall $\{x_3  =0 \}$. The flow is tangent to the inner wall
$\Gamma$, and the discontinuity of the tangent flow is given by the 
tangential gradient of $g$. Alternatively, if we wrote the analog of problem
 (\ref{BVP1}-\ref{Decay1}) in two dimensions, it could model a strike-slip fault in geophysics
where the equations of linear elasticity simplify to the scalar Laplacian.
In that case the scalar function $u$
models displacements in the direction orthogonal to a cross-section,
 $g$ models a slip, and normal derivatives model traction
\cite{ionescu2006inverse, ionescu2008earth, dascalu2000fault}. 
 A stability analysis for the inverse problem related to the two dimensional analog of 
(\ref{BVP1}-\ref{Decay1})  will be the subject of another study.  

\subsection{Existence and uniqueness for solutions to the direct problem}
  We seek solutions to problem (\ref{BVP1}-\ref{Decay1}) that have local $H^1$ regularity.  Due to the trace
	theorem (which is also valid in Lispchitz domains, \cite{ding1996proof, gagliardo1957caratterizzazioni}),
	we require the forcing term $g$ to be in
	  $\tilde{H}^{\f12} (\Gamma)$ which is the space of functions
	in $H^{\f12}(\p D)$ supported in $\ov{\Gamma}$.
As to $u$, we take it to be in the space 
\bea
 {\cal V} = 
\{ v \in H^1_{loc} (\RR^{3-}\sm \ov{\Gamma}): \nabla v,
 \f{v}{\s{ 1 + |\bx|^2}} \in
   L^2(\RR^{3-}\sm \ov{\Gamma})\}.
\eea
The natural norm $\|\nabla v \|_{L^2(\RR^{3-}\sm \ov{\Gamma})} +
 \| \f{v}{\s{ 1 + |\bx|^2}}\|_{L^2(\RR^{3-}\sm \ov{\Gamma})}$ on ${\cal V}$ is equivalent to the norm
 given by $\|\nabla v \|_{L^2(\RR^{3-}\sm \ov{\Gamma})} $ according to
Theorem A.1 in \cite{volkov2017reconstruction}.

\begin{thm}
    \label{UniquenessResult}
		Let $g$ be in $ \tilde{H}^{1/2}(\Gamma)$.
    Problem (\ref{BVP1}-\ref{BVP4}) has a unique solution $u$ in $\mathcal{V}$. 
		This solution satisfies the decay condition \eqref{Decay1}.
\end{thm}
\textbf{Proof:}
We first show uniqueness. Assume that $g=0$. 
Let $v$ be a compactly supported element in ${\cal V}$ such that $[v] =0$ across
$\Gamma$.
By \eqref{BVP2}, $\int_{\{x_3=0\}} (\p_{x_3} u)  v =0$.
Then, using Green's theorem and \eqref{BVP1}, \eqref{BVP3}, and  \eqref{BVP4},  we find that
$\int_{\RR^3-} \nabla u \cdot \nabla v =0$. 
By a density argument, we then infer that $\int_{\RR^3-} |\nabla u|^2=0$. 
But the definition of ${\cal V}$ precludes constant functions 
to be different from zero.
Thus $u=0$ and uniqueness is shown. \\
Existence can be shown by setting up a variational problem for $u$ in ${\cal V}$.
Instead, we choose 
to express $u$ as an integral over $\Gamma$ against an adequate Green function since this formulation will be crucial
further in this study. Denote,
\bean \label{fund}
\Phi(\bx,\by) = \f{1}{4 \pi} \f{1}{|\bx - \by|},
\eean
 the free space Green function
for the Laplacian. Then we can argue that
\bean \label{freeint}
 u'(\bx) = \int_{\Gamma} \nabla_{\by} \Phi (\bx,\by)  \cdot \bn(\by) g(\by) d \sigma (\by),
\eean
satisfies \eqref{BVP1}, \eqref{BVP3}, \eqref{BVP4}, \eqref{Decay1}, and is in 
${\cal V}$.
Indeed, as $g$ is in $ \tilde{H}^{1/2}(\Gamma)$, we can extend $g$ by zero to a function
in $ H^{1/2}(\p D)$ so (\ref{freeint}) equals
\bean \label{freeint2}
 u'(\bx) =  \int_{\p D} \nabla_{\by} \Phi (\bx,\by)  \cdot \bn(\by) g(\by) d \sigma (\by).
\eean
Even though $\p D$ is only Lipschitz regular, by Theorem 1 in \cite{costabel1988boundary},
$u'$ is a function in $H^1_{loc} (\RR^{3} \sm \ov{\Gamma})$
and by Lemma 4.1 in \cite{costabel1988boundary} the jump $[ u']$ across
$\Gamma$ is equal to $g$ almost everywhere, while the jump $[ \f{\p u' }{\p \bn}]$ is zero.
To find a solution to the PDE (\ref{BVP1}-\ref{BVP4})
we then set
\bean \label{intform}
 u(\bx) =  \int_{\Gamma} H (\bx,\by)   g(\by) d \sigma (\by),
\eean
where 
\bean \label{Hformula}
H (\bx,\by) = \nabla_{\by} \Phi (\bx,\by)  \cdot \bn(\by)  +
\nabla_{\by} \Phi (\ov{\bx},\by)  \cdot \bn(\by) ,
\eean
and $\ov{\bx} = (x_1, x_2, - x_3)$.
Then conditions (\ref{BVP1}-\ref{BVP4}) are satisfied,  $u(\bx) = O(\f{1}{|\bx|^2})$,
 and $\nabla u(\bx) = O(\f{1}{|\bx|^3})$, uniformly in 
$\f{\bx}{|\bx|}$, so $u$ is in ${\cal V}$. 
$\square$ 


\subsection{The crack inverse problem: formulation and uniqueness of solutions}
We now  prove a theorem stating that the crack inverse problem related to
problem (\ref{BVP1}-\ref{Decay1}) has at most one solution.
The data for the inverse problem is overdetermined boundary conditions over a portion
of the top plane $\{ x_3 = 0 \}$. Both the discontinuity $g$ and the crack $\Gamma$ are unknown 
in the inverse problem.  
Our theorem is valid for surfaces $\Gamma$
that are graphs of functions $\psi:R \ri \RR$.
We assume that $R $ is a bounded subset of $\RR^2$ such that 
it is equal to the interior of its closure and that the closure of $\Gamma$ is in
		 $\doubleR^{3-}$.
If we only assume that $\psi $ is Lipschitz regular, the crack inverse problem
may not be uniquely solvable. We show a counter-example in Appendix \ref{counter-example}.

\begin{thm}
\label{InverseProblemResult}
For $i=1,2$, let $R_i$ be an open and bounded subset of $\RR^2$ such that 
$R_i$ is equal to the interior of its closure.
    Let $\Gamma_i$   be the  surface
			defined by the graph of a  Lipschitz function $\psi_i:R_i \ri \RR$.
			Assume that the closure of $\Gamma_i$ is in
		 $\doubleR^{3-}$. 
		For  $i=1,2$, let 
 $u^i$ be the unique solution to  (\ref{BVP1}-\ref{Decay1}) with $\Gamma_i$ in place of $\Gamma$ and the jumps $g^i$ in  $\tilde{H}^{1/2}(\Gamma_i)$ in place of $g$ such that  
	$g^i$ has full support in  $\ov{\Gamma_i}$. 
	Let $V$ be a non-empty relatively open subset of $\{x_3=0\}$. 
	Assume that one of the three conditions hold,
	\begin{enumerate}[label=\arabic*., wide=0.5em,  leftmargin=*]
  \item $\psi_i$ is real analytic for $i=1,2$,
	\item $g^i$  is not equal to a constant on any relatively open subset of $\Gamma_i$,
	\item $\RR^{3-} \setminus \ov{\Gamma_1 \cup \Gamma_2}$ has only one connected component.
	\end{enumerate}
	If $u^1=u^2$ on $V$, then 
	$\Gamma_1=\Gamma_2$ and $g^1=g^2$ almost everywhere.	
\end{thm}
\textbf{Remark}: In practice, assumption 3  may be the most interesting. For example it is achieved if 
$\Gamma_i$ is planar as assumed further in this paper, or if  it is the union of two polygons such as in  the numerical reconstruction  
considered in  \cite{volkov2020stochastic}.

\textbf{Proof of Theorem \ref{InverseProblemResult}}:
Let $U=\RR^{3-} \setminus \ov{\Gamma_1 \cup \Gamma_2}$ and set $u= u^1 - u^2$ in $U$.
Given that 
$\Gamma_1$ and $\Gamma_2$ are bounded, 
$U$ has only one unbounded connected component which will be denoted by $U_1$.
Note that its boundary $\p U_1$ contains the top boundary 
$\{x_3 = 0\}$ since $\Gamma_1$ and $\Gamma_2$ are bounded away from that plane.
 Now, $u$ satisfies $\Delta u =0 $ in $U_1$ and $u=\p_{x_3} u =0$ on $V$,
so by the Cauchy Kowalevski theorem  and analytic continuation, $u$ is zero in $U_1$.
Let $U_j$, $j \neq 1 $, be a bounded connected component of $U$. 
Note that $\p U_j$ is not necessarily Lipschitz regular.
However, $u$ is in $H^1(U_j)$ and satisfies $\Delta u = 0$ in $U_j$ by construction.
We show in Appendix \ref{geom} that although $\p U_j$ may not be globally Lipschitz,
it is the union of two Lipschitz pieces, defined by $x_3= \psi_1(x_1, x_2)$ and
$x_3= \psi_2(x_1, x_2)$ where $(x_1,x_2)$ is in the closure of  a bounded and 
connected open subset $\Omega_j$ of $\RR^2$. 
We have $\psi_1(x_1, x_2)< \psi_2(x_1, x_2)$ for all $(x_1, x_2)$ in $\Omega_j$ or
$\psi_1(x_1, x_2)> \psi_2(x_1, x_2)$ for all $(x_1, x_2)$ in $\Omega_j$
as shown in Appendix  \ref{geom}, Lemma \ref{pos or neg}. 
Extending $\psi_1$ and $\psi_2$ 
to Lipschitz continuous functions on $\ov{R_1} $ and $\ov{R_2}$, 
 $\psi_1(x_1, x_2)= \psi_2(x_1, x_2)$, for $(x_1, x_2)$ on the boundary of $\Omega_j$
as shown in Appendix \ref{geom}, Lemma \ref{omega j}. 
 Due to relation  \eqref{BVP3}, 
$\aaa \frac{\p u^1}{\p \bn} \bbb$ is zero
on  the portion of $\p U_j$ where $x_3= \psi_1(x_1, x_2)$.
$\aaa \frac{\p u^2}{\p \bn} \bbb$ is also zero
on  that portion by regularity, so $\aaa \frac{\p u}{\p \bn} \bbb$ is zero.
A similar argument shows that $\aaa \frac{\p u}{\p \bn} \bbb$ is also zero on 
 the portion of $\p U_j$ where $x_3= \psi_2(x_1, x_2)$.
As $u$ was proven to be  zero in $U_1$ and the boundary
of $U_j$ is included in the boundary of $U_1$, we find that
  $ \frac{\p u}{\p \bn} $ is zero on $\p U_j$.
Green's theorem can now be applied to $u$ in $U_j$ to find that 
$\nabla u = 0$ in $U_j$. As $U_j$ is connected, 
$u$ is equal to some constant $C_j$ in $U_j$. If $j=1$, we set $C_j=0$. \\
Let $(y_1,y_2)$ be in $R_1$. Arguing by contradiction, assume that $(y_1,y_2) \notin \ov{R_2}$.
Then there is an open ball $B((y_1,y_2),\alpha)$ in $\RR^2$ such that  $B((y_1,y_2),\alpha) \subset R_1$ and
$B((y_1,y_2),\alpha) \cap R_2 = \emptyset$. Let $\by = (y_1,y_2, \psi_1(y_1,y_2))$.
Then there is an open ball $B(\by,\beta)$ in $\RR^3$ such that  $B(\by,\beta) \cap \Gamma_2 = \emptyset$.
Clearly $B(\by,\beta) \setminus \ov{\Gamma_1}  \subset U_1$, thus $[u] =0$ in $B(\by,\beta) \cap \Gamma_1 $. 
It is also clear that $[u^2] =0$ in $B(\by,\beta) \cap \Gamma_1 $, 
so we infer that $[u^1] =0$ in $B(\by,\beta) \cap \Gamma_1 $:
this contradicts the assumption that $[u^1]$ has full support in $ \Gamma_1$. We have thus shown 
that $R_1 \subset \ov{R_2}$. By switching the roles of $R_1$ and $R_2$, this shows that 
 $\ov{R_1} = \ov{R_2}$. Due to our assumption on $R_1$ and $R_2$, this implies that $R_1 = R_2$.\\
Now, let $R_3$ be the interior of the set 
$\{(x_1,x_2) \in R_1: \psi_1(x_1,x_2) = \psi_2(x_1,x_2)\}$. 
Since the half-line 
$\{(x_1,x_2, t  \psi_1(x_1,x_2) ): t \in (1,\infty)\}$
and the line segment $\{(x_1,x_2, t  \psi_1(x_1,x_2) ): t \in (0,1)\}$
 are  in $U_1$ for all $(x_1, x_2)$ in $R_3$,
 we can argue that  $[u] = 0$ in $R_3 \times \psi_1(R_3) $, and thus  $g^1=g^2$ in this set. In particular, if $\Gamma_1 = \Gamma_2$
then $g^1=g^2$ almost everywhere. \\
There remains to show that $R_3=R_1$.
Fix  $\by = (y_1,y_2, \psi_1(y_1,y_2))$ for $(y_1,y_2)$ in $R_1$. 
Arguing by contradiction,
assume that 
$\psi_1(y_1,y_2) \neq \psi_2(y_1,y_2)$. Then, by continuity, $\psi_1(z_1,z_2) \neq \psi_2(z_1,z_2)$, for all 
$(z_1,z_2) $ in an open neighborhood $W\subset R_1$ of $(y_1,y_2)$. 
Since $u^2$ is regular
 near $\by$, we can argue that $[u^2(y_1,y_2, \psi_1(y_1,y_2))]=0$ 
for  almost all $(y_1,y_2)$  in 
$W$,  thus 
$[u^1 (y_1,y_2, \psi_1(y_1,y_2))] = \pm C_j$ for some $j$ and  almost all $(y_1,y_2)$  in 
$ W$.
If condition 2 holds, this is a contradiction. If condition 3 holds, then  $C_j=C_1=0$, and this contradicts that
$g^1 $ has full support in $\Gamma_1$. If condition 1 holds, given that there are no piecewise constant functions
in $\tilde{H}^{\f12} (\Gamma_i)$ other than zero, if $C_j \neq 0$ for some $j$, 
  there is an open set $S$ of $\RR^{3-}$  such that 
	$\Gamma_1 \cap S \neq \emptyset$ and
	$[u^1]$ is not constant in $\Gamma_1 \cap S$. But given what we showed earlier, 
	we must have that $\Gamma_1 \cap S = \Gamma_2 \cap S$,
	that is, the functions  $\psi_1$ and $\psi_2$ coincide on the open subset of $\RR^2$
	$\{ (x_1,x_2) \in \RR^2: (x_1,x_2, x_3) \in \Gamma_1 \cap S \}$.
	Thus if condition 1 holds, 
	$\psi_1=\psi_2$ everywhere in $R_1=R_2$.  $\square$ \\
	
\subsection{High order jump formulas for derivatives of the  double layer potential}
To prove  stability results for the crack inverse problem, we will argue by contradiction.
From the assumption that the  crack inverse problem is not 
Lipschitz continuous, we will 
form a PDE on $\Gamma$ for the slip $g$, which we will show to  only have the trivial solution, 
pointing to a contradiction. 
To form that PDE for $g$, we will need jump formulas for the integral expression 
\eqref{intform} for $\bx$ across $\Gamma$. 
Although the  formula for the jump of \eqref{intform} across $\Gamma$ is well known, 
formulas for the jumps of its first and second derivatives are less known, so it is worth stating them in detail. A proof is provided in Appendix \ref{appen prove jump}.  

With $\Phi$ as in \eqref{fund} and $\bv$ a vector in $\RR^3$, we set
for $\bx, \by$ in $\RR^{3-}$ such that  $\bx \neq \by$,
\bean
 G(\bx,\by,\bv) = 
\nabla_\by \Phi (\bx,\by) \cdot \bv, \label{Gdef} \\
H(\bx,\by,\bv) = 
\nabla_\by \Phi (\bx,\by) \cdot \bv,
+\nabla_\by \Phi (\ov{\bx},\by) \cdot \bv,
\label{Hdef}
\eean
where as previously $\ov{\bx} = (x_1,x_2, -x_3)$.\\
For an open subset $V$ of $\RR^d$, we denote $C^\infty_c(V)$
the space of smooth functions compactly supported in $V$. 
If $\Gamma$ is an open planar surface in $\RR^3$, it can be mapped by
rotation and translation to an open subset of $\RR^2$ and this mapping can be used to define
$C^\infty_c(\Gamma)$ by pullback. 

\begin{lem} \label{thirdjumplem}
    Let $\Gamma$ be a planar open surface in $\doubleR^3$ with surface element
		$dS$, and unit normal vector $\bn$. Let 
	 $\bt$ be a fixed unit vector parallel to $\Gamma$.  Let 
		$g$ be in $L^2(\Gamma)$ and  $\phi$ in $C^\infty_c(\Gamma)$. 
		For $\epsilon > 0$, $\bx$, $\by$ in $\Gamma$, and $\bp$ in $\RR^3$ define the notation,
		\bea
		G^{+-}(\bx,\by,\bp) = G(\bx + \epsilon \bn,\by,\bp) - G(\bx - \epsilon \bn,\by,\bp).
		\eea
		The following jump formulas across $\Gamma$ hold, 
\begin{equation}
    \label{jumpformula31}
   \lim_{\epsilon \ri 0} \int_\Gamma \int_\Gamma G^{+-}(\bx,\by,\bn) g(y)
		dS(\by) \phi (\bx) dS (\bx) =\int_\Gamma g(\bx)  \phi (\bx)  dS (\bx), 
\end{equation}
\begin{equation}
    \label{jumpformula32}
  \lim_{\epsilon \ri 0}   \int_\Gamma \int_\Gamma \frac{\partial_{\by}G^{+-}}{\p \bt}(\bx,\by,\bn)g(\by)
	 dS(\by) 	\phi (\bx) dS (\bx)=
		\int_\Gamma g(\bx) \p_\bt \phi(\bx) dS (\bx), 
\end{equation}
\begin{equation}
    \label{jumpformula33}
   \lim_{\epsilon \ri 0}  \int_\Gamma \int_\Gamma G^{+-}(\bx,\by,\bt)g(\by) dS(\by) \phi (\bx) dS (\bx)=0,
\end{equation}

\begin{equation}
    \label{jumpformula34}
   \lim_{\epsilon \ri 0}  \int_\Gamma \int_\Gamma \f{\p_\by G^{+-}}{\p \bn}(\bx,\by,\bn)g(\by)
		dS(\by)  \phi (\bx) dS (\bx) =0.
\end{equation}

\begin{equation}
    \label{jumpformula342}
  \lim_{\epsilon \ri 0}   \int_\Gamma   \left(\f{\p_\bx}{\p \bn} \int_\Gamma G^{+-}(\bx,\by,\bn)g(\by)
		dS (\by)  \right) \phi (\bx) dS(\bx)=0,
\end{equation}

\begin{equation}
    \label{jumpformula35}
 \lim_{\epsilon \ri 0}    \int_\Gamma \left( \f{\p_\bx}{\p \bn}
		\int_\Gamma G^{+-}(\bx,\by,\bt)g(\by)dS (\by) \right)  \phi (\bx) dS(\bx)=
		-\int_\Gamma g(\bx) \partial_{\bt}\phi(\bx) dS (\bx),
\end{equation}

\begin{equation}
    \label{jumpformula36}
    \lim_{\epsilon \ri 0}   \int_\Gamma   \left(\f{\p_\bx}{\p \bn}\int_\Gamma \f{\p_\by G^{+-}}{\p \bn}
		(\bx,\by,\bn)g(\by) dS(\by) \right) \phi (\bx)dS (\bx) =
	\int_\Gamma  g(\bx) 	\Delta_\Gamma \phi(\bx) dS (\bx).
\end{equation}

\begin{equation}
    \label{jumpformula343}
    \lim_{\epsilon \ri 0} \int_\Gamma \left(\f{\p_\bx}{\p \bn}\int_\Gamma \f{\p_\by G^{+-}}{\p \bt}
		(\bx,\by,\bn)g(\by)dS(\by) \right) \phi (\bx)  dS (\bx)= 0.
\end{equation}

\end{lem}

\section{Lipschitz stability results for the half space crack inverse problem
in  case of  planar geometries}
For a  bounded relatively open  set $R$ in the plane $\{x_3=0 \}$  and
$\bm= (a,b,d)$ in $\RR^3$, 
we define 
the surface
 \bea
\Gamma_{\bm}=\{(x_1,x_2,ax_1+bx_2+d):(x_1,x_2)\in R \}.
\eea 
 Let $B$ be a set of $\bm$ 
in $\RR^3$
such that $\Gamma_{\bm}\subset\doubleR^{3-}$.
We assume that $B$ is  closed and bounded.  It follows  that
there is a positive constant $\beta$ such that
\begin{equation}
    \label{BProperty}
   \mbox{dist }(\Gamma_\bm,\{x_3=0\}) \geq \beta,
\end{equation}
for all $\bm$ in $B$, where $\mbox{dist }$ is the distance between sets.
We choose the unit normal vector on $\Gamma_\bm$ to be
 $\bn=\frac{(-a,-b,1)}{\sqrt{a^2+b^2+1}}$. The surface element on
$\Gamma_\bm$ will be
 denoted by
$\sigma d x_1 d x_2=\sqrt{a^2+b^2+1} d x_1 d x_2$. 
 Using \eqref{intform} and \eqref{Hdef}, we 
introduce the linear compact operator mapping jumps $g$ across 
$\Gamma_\bm $ to Dirichlet data on $V$,
\bean
 A_\bm:H_0^1(R)\to L^2(V), \no \\
      (A_\bm g )(\bx) =\int_R H(x_1,x_2,0,y_1,y_2,ay_1+by_2+d,\bn)
				g(y_1,y_2)\sigma
				dy_1dy_2,    \label{AmDef}
\eean
where $\bx=(x_1,x_2,0)$ is in $V$. 

We know from Theorem \ref{InverseProblemResult} that $A_\bm$ is injective.
 Fix a  non-zero $h$ in $H_0^1(R)$ and define the 
function $\phi:B\to L^2(V)$ by 
\begin{equation}
    \label{phiDef}
    \phi(\bm)= A_\bm h.
\end{equation}
Due to the regularity of the Green's function $H$ and 
(\ref{BProperty}), we know that $\phi$ is analytic in $\bm$. 
Theorem \ref{InverseProblemResult} implies that $\phi$ is injective. We will show that the inverse of $\phi$ defined on $\phi(B)$ and valued in $B$ is of class $C^1$ by applying
the  inverse function Theorem.
As a consequence, $\phi^{-1}$ is Lipschitz continuous, which will yield
our first stability estimate.
\begin{thm}
\label{InverseFunctionTheorem}
Fix a non-zero $h$ in $H_0^1(R)$ and define the function $\phi$ from $B$ to $L^2(V)$ by \eqref{phiDef}. Then there is a positive constant $C$ such that 
\begin{equation}
    \label{LipschitzEst}
    C|\bm- \bm'|\le \norm{\phi(\bm)-\phi(\bm')}_{L^2(V)}
\end{equation}
for all $\bm$ and $\bm'$ in $B$.
\end{thm}
\textbf{Proof:}
Let ${\cal B} $ be an open set of $\RR^3$ such that
$ B \subset {\cal B} $ and for all $\bm$ in ${\cal B} $ 
\begin{equation}
    \label{BProperty2}
   \mbox{dist }(\Gamma_\bm,\{x_3=0\}) > \f{\beta}{2},
\end{equation}
where $\beta$ is as in \eqref{BProperty}.
Arguing by contradiction, suppose that there is an $\bm$ in ${\cal B} $ 
 such that $\nabla \phi(\bm)$ does not have full rank. 
Then there is a non-zero vector $(\gamma_1,\gamma_2,\gamma_3)$ such that, 
\begin{equation}
    \label{LinearDep}
    \gamma_1\pfrac{}{a}\phi(\bm)+\gamma_2\pfrac{}{b}\phi(\bm)
		+\gamma_3\pfrac{}{d}\phi(\bm)=0.
\end{equation} We note that $ \bn \sigma$ simplifies to $(-a,-b,1)$. 
Since 
$H(\bx,\by,\bn)$
is linear in $\bn$, 
we can apply the chain rule with $y_3=ay_1+by_2+d$ on $\Gamma_\bm$ to find that
\begin{equation}
    \begin{split}
        \label{GreenDeriva}
        \pfrac{}{a}H(\bx,\by,\bn \sigma)&=
				\pfrac{y_3}{a}(\partial_{y_3}H)(\bx,\by,\bn\sigma)+
				H(\bx,\by,\pfrac{(\bn\sigma)}{a})\\
        &=y_1(\partial_{y_3}H)(\bx,\by,\bn \sigma)-H(\bx,\by,\bev_1).
    \end{split}
\end{equation}
Similarly,
\begin{equation}
    \begin{split}
        \label{GreenDerivb}
        \pfrac{}{b}H(\bx,\by,\bn\sigma)
        &=y_2(\partial_{y_3}H)(\bx,\by,\bn \sigma)-H(\bx,\by,\bev_2),
    \end{split}
\end{equation}
and
\begin{equation}
    \begin{split}
        \label{GreenDerivd}
        \pfrac{}{d}H(\bx,\by,\bn \sigma)
        &=(\partial_{y_3}H)(\bx,\by,\bn \sigma).
    \end{split}
\end{equation}
For $(y_1, y_2)$ in $R$, define $f(y_1,y_2)=\gamma_1y_1+\gamma_2y_2+\gamma_3$. Since $h$ is independent of $a$, $b$, and $d$, we substitute (\ref{GreenDeriva}-\ref{GreenDerivd}) into \eqref{LinearDep} to find that
\begin{equation}
    \begin{split}
        \label{LinearDep2}
        &\int_R (\partial_{y_3}H)(\bx,y_1,y_2,ay_1+by_2+d,\bn)h(y_1,y_2)f(y_1,y_2)\sigma dy_1dy_2\\
        &\qquad-\int_R H(\bx,y_1,y_2,ay_1+by_2+d,\nabla f)h(y_1,y_2)dy_1dy_2=0,      
    \end{split}
\end{equation}
for all $\bx$ in $V$. Set $w(\bx)$ to be the left hand side of \eqref{LinearDep2} 
where $\bx$ has been extended to $\doubleR^{3-}\sm\Gamma_\bm$. 
We will now show that $w$ is zero in $\doubleR^{3-}\sm\Gamma_\bm$. 
Since $\partial_{y_{3}}$ and $\partial_{x_i}$ commute, 
we know from the definition of the Green's function that $w$ satisfies the Laplace Equation in
$\doubleR^{3-}\sm\Gamma_\bm$.
 We also know that $w$ is $0$ on $V$ thanks to \eqref{LinearDep2}. 
By construction of the Green's tensor $H$, we know that for any $x$ on the plane $x_3=0$, any $\by\in \doubleR^{3-}$ and any fixed vector 
$\bp\in \doubleR^{3}$,
using the notation
\eqref{Hdef},
 \[\partial_{x_3}H(\bx,\by,\bp)=0.\]
Thus we can take a $\partial_{y_3}$ derivative and commute it with 
$\partial_{x_3}$ to obtain
\[\partial_{x_3}\partial_{y_3}H(\bx,\by,\bp)=0.\]
It follows that $\partial_{x_3}w$ is also zero in $V$. 
As $w$ is zero in $V$,
and $\Delta w =0 $ in $\RR^{3-} \sm \ov{\Gamma_\bm}$,
 $w=0$ everywhere in $\doubleR^{3-}\sm\Gamma_m$ 
by the Cauchy Kowalevski theorem. 
In particular, the jump of $w$ across $\Gamma_m$ must be zero. 
As mentioned  earlier, $H(\bx,\by,\bp)-G(\bx,\by,\bp)$ is smooth for
 any $\bx, \by$ in $\doubleR^{3-}$
and any fixed vector $\bp$ in $\doubleR^3$, therefore the jump 
across $\Gamma_\bm$ of 
   \bean
        \label{LinearDep3}
     \int_R (\partial_{y_3}G)(\bx,y_1,y_2,ay_1+by_2+d,\bn)
			h(y_1,y_2)f(y_1,y_2)\sigma dy_1dy_2  \no \\
     -\int_R G(x,y_1,y_2,ay_1+by_2+d,\nabla f)h(y_1,y_2)dy_1dy_2    
\eean
is also zero. We now consider two cases.\\
\textbf{Case 1:} $(a,b) \neq 0$.\\
In this case the fault $\Gamma_m$ is not horizontal.
Let $\alpha$ and $\beta$ be in $\doubleR$ and $\bt$
a unit vector in $\doubleR^3$ parallel to
  to $\Gamma_m$ such that,
 \bean \label{alphabeta}
\bev_3=\alpha\bn+\beta\bt.
\eean
 Then (\ref{LinearDep3}) can be rewritten as
\begin{equation}
\begin{split}
\label{corrected}
    &\alpha\int_R  \f{\p_\by G}{\p \bn}(\bx,y_1,y_2,ay_1+by_2+d,
		\bn)h(y_1,y_2)f(y_1,y_2)
		\sigma dy_1dy_2\\
    &\qquad+\beta
		\int_R \f{\p_\by G}{\p \bt} (\bx,y_1,y_2,ay_1+by_2+d,\bn )
		h(y_1,y_2)f(y_1,y_2)\sigma dy_1dy_2\\
		 &\qquad - \int_R G(\bx,y_1,y_2,ay_1+by_2+d,\nabla f)h(y_1,y_2)dy_1dy_2,
\end{split}
\end{equation}
which  is again zero for $\bx$ in  
$\doubleR^{3-}\setminus \overline{\Gamma_\bm}$.
We now apply lemma \ref{thirdjumplem}: to be more precise,  respectively 
formulas (\ref{jumpformula31}) to (\ref{jumpformula34})
to obtain the weak equation
\begin{equation}
    \label{eqonh}
- \beta \frac{\partial}{\partial \bt} (h f)
- \frac{1}{\sigma} (\nabla f \cdot \bn ) h =0 .
\end{equation}
Since in case 1 $(a,b) \neq 0$, $\bev_3 \neq \bn$, so we infer
that $\beta \neq 0$. 
As $h$ is in $H_0^1(R)$, equation (\ref{eqonh}) implies that 
$h$ is zero: this is due to lemma 3.3 in \cite{triki2019stability}.
Thus we have contradicted the assumption $h \neq 0$. \\
\textbf{Case 2:} $a$ and $b$ are zero.\\
In that case the fault $\Gamma_\bm$ is  horizontal. Recalling that $w(\bx)=0$ 
for all
$x$ in $\RR^{3-}\sm{\ov{\Gamma_\bm}}$, its $x_3$ derivative is also zero in the same open set thus the jump
across $\Gamma_m$ of 
\bean
        \label{d3w}
      \p_{x_3}w(\bx)  & = &\p_{x_3} \int_R (\partial_{y_3}G)(\bx,y_1,y_2,d,
			\bev_3)
			h(y_1,y_2)f(y_1,y_2)dy_1dy_2  \no \\
      &  &\qquad-\p_{x_3} \int_R G(\bx,y_1,y_2,d,\nabla f)h(y_1,y_2)dy_1dy_2    
\eean
is zero. We then apply the jump formulas  (\ref{jumpformula42w}-\ref{jumpformula6w}) 
from  Lemma
\ref{secondjumplem} to obtain the weak equation,
$
 \Delta (f  h ) - \nabla f \cdot \nabla h = 0
$.
Since $\Delta f =0$, this simplifies to
$  
\mbox{div } (f \nabla h )= 0,
$
in other words
\bean  \label{ellipeq}
\int_R f \nabla h \cdot \nabla \varphi= 0,
\eean
for all $\varphi$ in $H^1_0(R)$. Recall that $h$ is in $H^1_0(R)$ and $f(x_1, x_2) =
\gamma_1 x_1 + \gamma_2 x_2 + \gamma_3 $ with
$(\gamma_1 ,  \gamma_2 , \gamma_3 ) \neq 0 $.
If $f$ does not change signs in $R$, \eqref{ellipeq} clearly implies that $h$ is zero.
If $f$  changes signs in $R$, then there is an open set $R^+$ in $R$ where $f>0$, 
and an open set $R^-$ where $f<0$. $R^+$ and $R^-$ are separated by a line.
Without loss of generality, we can change coordinate systems  in $\RR^2$
by rotation and translation
 so that this line becomes the line with equation $x_1=0$ and we can rescale
$f$ so that it becomes $f(x_1,x_2) =x_1$.  $R^+$ is  then  the set
$\{ (x_1, x_2) \in R: x_1 >0 \}$.  For $\epsilon >0$, let
$p_\epsilon$ be the function from $\RR^2$ to
$[0,1]$ be defined by
\bea
p_\epsilon (x_1, x_2)= \left\{
\begin{array}{ll}
1 & \mbox{ if } x_1 \geq \epsilon, \\
\ds \f{x_1}{\epsilon} & \mbox{ if } 0 \leq  x_1 < \epsilon, \\
0 & \mbox{ if }  x_1 <0.
\end{array}
\right.
\eea
$p_\epsilon$ is Lipschitz continuous and $p_\epsilon h$ is in $H^1_0(R)$.
By \eqref{ellipeq}, 
$
\int_R f \nabla h \cdot \nabla (p_\epsilon h)= 0
$, which we write as 
\bean \label{two terms}
\int_R p_\epsilon f |\nabla h|^2   \,
+ \,
\int_R   h \nabla h \cdot f\nabla p_\epsilon = 0.
\eean
As $\epsilon$ tends to zero the limit of the first term in \eqref{two terms}
is $\int_R f^+ |\nabla h|^2$, where $f^+$ is the positive part of $f$.
To tackle the second term,  we note that $h$ is also in $L^3(R)$ by the Sobolev embedding theorem,
thus by the  generalized H\"older's  inequality, 
$h \nabla h$ is in $L^{\f{6}{5}}(R)$, as $\f12 + \f13 = \f56$.
It follows, using  H\"older's  inequality, that 
\bea
|\int_R   h \nabla h \cdot f\nabla p_\epsilon |
\leq \| h \nabla h \|_{L^{\f{6}{5}}(R)} (\int_R   | f\nabla p_\epsilon|^6)^{\f16}.
\eea
But now
\bea
 \int_R   | f\nabla p_\epsilon|^6 = \int_{R \cap \{(x_1, x_2): 0 < x_1 <\epsilon\}} 
(\f{x_1}{\epsilon})^6 \, dx_1 dx_2  
\eea
is $O(\epsilon)$ thanks to the rescaling $y_1= \f{x_1}{\epsilon}$, thus we have shown that the second term in \eqref{two terms} converges to zero, thus
$\int_R f^+ |\nabla h|^2 =0$.  Similarly,  we can show that 
$\int_R f^- |\nabla h|^2 =0$, so altogether, $
\int_R  |f| |\nabla h|^2 =0$.
As $f$ is a non-zero affine function, this implies that $h$
is zero almost everywhere: contradiction. \\
In summary, we have proved that $\nabla\phi(\bm)$ has full rank for any $\bm$ in
${\cal B}$.
The inverse function theorem guarantees that $\phi$ defines a $C^1$ diffeomorphism from an open neighborhood $U_\bm$ of $\bm$ in ${\cal B}$ to its image by $\phi$ in $L^2(V)$. 
Arguing by contradiction, assume that estimate \eqref{LipschitzEst} does not hold. 
Then there are two sequences $\bp_n$ and $\bq_n$ 
in $B$
such that $\bp_n\ne \bq_n$ and
\begin{equation}
    \label{ProblemEq}
    \lim_{n \ri\infty}\frac{\norm{\phi(\bp_n)-\phi(\bq_n)}_{L^2(V)}}{\a{\bp_n-\bq_n}}=0.
\end{equation}
As $B$ is compact, 
without loss of generality we may assume that 
$\bp_n$ converges to some $\tilde{\bm}$ in $B$ and
$\bq_n$ converges to some $\tilde{\tilde{\bm}}$ in $B$.
If $\tilde{\bm} \neq \tilde{\tilde{\bm}}$, this contradicts 
the uniqueness Theorem \ref{InverseProblemResult}, so we claim that 
both $\bp_n$ and $\bq_n$ converge to $\tilde{\bm}$ in $B$.
But we proved that $\phi$ defines a $C^1$ diffeomorphism
in a open neighborhood $U_{\tilde{\bm}}$ of $\tilde{\bm}$, thus there exists 
a positive constant $C_{\tilde{\bm}}$ such that for all $n$ large enough
\bea
|\bp_n - \bq_n| \leq C_{\tilde{\bm} }
\| \phi(\bp_n)-\phi(\bq_n) \|_{L^2(V)} ,
\eea
which contradicts \eqref{ProblemEq}.
$\square$

Stability Theorem \ref{InverseFunctionTheorem}
holds for a fixed jump $h$ in $H_0^1(R)$,  but
jumps are unknown in the crack inverse problem. The following result
addresses this important issue. A Lipschitz stability statement in $h$ is not possible  
due the compactness of the operator $A_{\bm}$, however even
if the jump $h$ is variable, a Lipschitz stability result in the geometry parameter 
$\bm$ holds as shown below.   
\begin{thm}\label{stability2}
Let $B$ be a closed and bounded subset of geometry  parameters $\bm$ in $\RR^3$ such that 
  the distance condition \eqref{BProperty} holds. Let $A_\bm$ be the crack to boundary 
	operator defined by \eqref{AmDef}.
Fix a non-zero $h_0$ in $H^1_0(R)$ and $\bm_0$ in $B$.
There exists a positive constant $C$ such that
\bean \label{inf formula}
\inf_{h \in H^1_0(R)} \| A_\bm h -  A_{\bm_0} h_0  \|_{L^2(V)}  \geq
 C |\bm - \bm_0|,
\eean
for all $m$ in $B$.
\end{thm}
\textbf{Proof:} 
Let $P_\bm$ be the orthogonal projection onto $\ov{{\cal R}(A_\bm)}$, the closure of the
 range 
of $A_\bm$  in $L^2(V)$ 
for $\bm$ in $B$. 
We first claim that,
\bean \label{proj_est}
\| P_\bm - P_{\bm_0} \| = O(|\bm - \bm_0|).
\eean
Indeed, since the nullspaces of $A_\bm^*$ and 
$A_\bm A_\bm^*$ are equal, taking their orthogonal we find that
 $\ov{{\cal R}(A_\bm)}= \ov{{\cal R}(A_\bm A_\bm^*)}$. 
The orthogonal projection on 
$\ov{{\cal R}(A_\bm A_\bm^*)}$
 can be expressed by the contour integral  \cite{kato2013perturbation}
\bea
P_\bm = \f{1}{2 i\pi} \int_{{\cal C}} (   \zeta I - A_\bm A^*_\bm)^{-1} d \zeta,
\eea
where ${\cal C}$ is the circle in the complex plane centered at the origin with radius
$\|A_{\bm} A^*_{\bm} \| +1$.
Since 
$\| A_\bm A^*_\bm - A_{\bm_0} A^*_{\bm_0}\| = O(|\bm-\bm_0|)$,
 estimate \eqref{proj_est}
 holds.\\
Arguing by contradiction, assume that there is a sequence $\bm_n$ in $B$ 
 and a sequence $h_n$ in $H^1_0(R)$ such that 
$\bm_n \neq \bm_0$ for all $n \geq 1$ and
\bean \label{little}
\| A_{\bm_0} h_0 - A_{\bm_n} h_{m_n}   \|_{L^2(V)} = o (|\bm_n - \bm_0|).
\eean
Since $I - P_{\bm_n}$ is an orthogonal projection, 
\bea
\| A_{\bm_0} h_0 - A_{\bm_n} h_{\bm_n}   \|_{L^2(V)}  \geq 
\| (I - P_{\bm_n} ) (A_{\bm_0} h_0 - A_{\bm_n} h_{\bm_n} )  \|_{L^2(V)}  ,
\eea
thus by  \eqref{little},
\bea
\| (I - P_{\bm_n} ) A_{\bm_0}  h_0 \|_{L^2(V)}  = o (|\bm_n - \bm_0|),
\eea
or equivalently,
\bean \label{equiv}
\| (I - P_{\bm_n} ) ( A_{\bm_n} -A_{\bm_0} ) h_0 \|_{L^2(V)}  
= o (|\bm_n - \bm_0|).
\eean
By \eqref{proj_est},
\bean \label{double}
\| ( P_{\bm_0}- P_{\bm_n} ) ( A_{\bm_n} -A_{\bm_0} ) h_0 \|_{L^2(V)}  
= o (|\bm_n - \bm_0|),
\eean
thus combining \eqref{equiv} and \eqref{double}, we find,
\bean \label{wefind}
\| (I - P_{\bm_0} ) ( A_{\bm_n} -A_{\bm_0} ) h_0 \|_{L^2(V)}  
= o (|\bm_n - \bm_0|).
\eean
Equivalently,
\bean \label{proj}
 (I - P_{\bm_0})  \f{(A_{\bm_n} - A_{\bm_0})}{|\bm_n - \bm_0|} h_0  =
 o (1).
\eean
As $\f{\bm_n - \bm_0}{|\bm_n - \bm_0|}$ is a sequence on the unit sphere 
of $\RR^3$, after
possibly extracting a subsequence we may assume that it converges to some
$\bq$ with $|\bq|=1$.
To  keep notations consistent with those used in the  
proof of Theorem \eqref{InverseFunctionTheorem},
we set $q=(\gamma_1,\gamma_2, \gamma_3)$,
 $\bm_0 = (a,b, d)$,
and $f(y_1,y_2)=\gamma_1y_1+\gamma_2y_2+\gamma_3$.
We define the linear  operator 
\bean \label{partial der def}
 \p_{\bq} A_{\bm_0} :H_0^1(R)\to L^2(V) \no \\
     \p_{\bq} A_{\bm_0} = \gamma_1\f{\p}{\p a} A_{\bm_0} +
			\gamma_2\f{\p}{\p b} A_{\bm_0}+\gamma_3\f{\p}{\p d} A_{\bm_0}.
\eean
Taking the limit as $n \ri \infty$ in (\ref{proj}) we find,
\bean \label{preg0}
 (I - P_{\bm_0}) \p_{\bq} A_{\bm_0}  h_0 = 0,
\eean
thus, there is a $g_0$ in $H^1_0(R)$ such that,
\bean \label{with rhs}
\p_{\bq} A_{\bm_0}   h_0 -  A_{m_0}  g_0 =0.
\eean
We now show that \eqref{with rhs} implies that 
$ h_0 = 0$, which is a contradiction.
Based on \eqref{with rhs},  we define
    \bean
     w'(\bx)& =& \int_R (\partial_{y_3}H)(\bx,y_1,y_2,ay_1+by_2+d,\bn)h_0(y_1,y_2)f(y_1,y_2)\sigma dy_1dy_2 \no\\
      &&  -\int_R H(\bx,y_1,y_2,ay_1+by_2+d,\nabla f)
			h_0(y_1,y_2)dy_1dy_2 \no\\        
	&&			- \int_R H(\bx,y_1,y_2,ay_1+by_2+d,\bn)g_0(y_1,y_2)\sigma dy_1dy_2,
	\label{w'def}
    \eean
for $\bx$ in $\RR^{3-} \sm \ov{\Gamma_{\bm_0}}$, 
$w'(\bx) $ is zero for $\bx$ in $V$ because of relation \eqref{with rhs}.
  As $\p_{x_3}w'(\bx) $ is also
 zero for $\bx$ in $V$, by the same argument as in the proof of Theorem 
\ref{InverseFunctionTheorem}, we infer 
that  $w'(\bx) =0 $ for all  $\bx $ in $\RR^{3-} \sm \ov{\Gamma_{\bm_0}}$.
Now, unlike, in the proof of  Theorem \ref{InverseFunctionTheorem}, we have
to contend with the inhomogeneity presented by the $g_0$ term in 	\eqref{w'def}. 
To eliminate that term, we take the jump $ [\p_{\bn} w'] $
across $\Gamma_{\bm_0}$. Once again we use that $H-G$ is smooth
in $\RR^{3-} \times \RR^{3-}$, thus we only need to consider the jump of
the normal derivative across $\Gamma_{\bm_0}$ of 
\bea \ds
    \alpha\int_R  \f{\p_\by G}{\p \bn}(\bx,y_1,y_2,ay_1+by_2+d,
		\bn)h_0(y_1,y_2)f(y_1,y_2)
		\sigma dy_1dy_2\\
   +\beta
		\int_R \f{\p_\by G}{\p \bt} (\bx,y_1,y_2,ay_1+by_2+d,\bn )
		h_0(y_1,y_2)f(y_1,y_2)\sigma dy_1dy_2\\
		- \int_R G(\bx,y_1,y_2,ay_1+by_2+d,\nabla f)h_0(y_1,y_2)dy_1dy_2 \\
	- \int_R G(\bx,y_1,y_2,ay_1+by_2+d,\bn)g_0(y_1,y_2)\sigma dy_1dy_2,
 \eea
where $\alpha, \beta, \bt $ are as in \eqref{alphabeta}.
To compute this  jump, we respectively apply  formulas 
\eqref{jumpformula36}, \eqref{jumpformula343},
\eqref{jumpformula342} and  \eqref{jumpformula35},
and \eqref{jumpformula342}. Doing so, we obtain
the weak  equation on $\Gamma_{\bm_0}$,
\bean
\label{eqm0}
\alpha \Delta_{\Gamma_{\bm_0}} (\tilde{f}  \tilde{h_0}) - \f{1}{\sigma} 
\nabla_{\Gamma_{\bm_0}} \tilde{f}   \cdot \nabla_{\Gamma_{\bm_0}} \tilde{h_0}  = 0,
\eean
where $\tilde{h_0}$ is the function in $H_0^1(\Gamma_{\bm_0} )$
defined by $\tilde{h_0} (x_1,x_2, ax_1 + b x_2 +d) = h_0(x_1, x_2)$ and 
$\tilde{f}$ is
defined from $f$ likewise.
From \eqref{alphabeta},
\bea
\alpha= \bn \cdot \bev_3 = \f{1}{\s{a^2+b^2+1}},
\eea
thus $\alpha = \f{1}{\sigma}$. 
Since $ \Delta_{\Gamma_{\bm_0}} \tilde{f} =0 $,
\eqref{eqm0} is equivalent to writing that
$\int_{\Gamma_{\bm_0}} \tilde{f}
\nabla_{\Gamma_{\bm_0}}  \tilde{h_0}  \cdot 
\nabla_{\Gamma_{\bm_0}} \varphi = 0 $, for all $\varphi$ in $H_0^1(\Gamma_{\bm_0} )$.
Just as in the case of variational problem \eqref{ellipeq}, we can argue from there 
that $\tilde{h_0} $ must be zero: 
this contradicts our assumption on $h_0$ thus 
 \eqref{little} cannot hold.\\
Therefore, for every $\bm_0$ in $B$, there is a positive constant
$C_{\bm_0}$ such that for all $\bm$ in $B$ and all $\bh$ in $H^1_0(R)$,
\bean \label{withCm0}
\| A_{\bm_0} h_0 - A_{\bm} h_{m}   \|_{L^2(V)} 
\geq C_{\bm_0} |\bm_0- \bm|.
\eean 
$\square$ \\

\vskip 10pt
Can the constant $C_{\bm_0}$ in (\ref{withCm0}) be chosen 
uniformly in $h_0$ and $\bm_0$, 
or at least uniformly  for $h_0$ in a bounded set?
The following theorem addresses this question.
Since the injective operator $A_{\bm_0}$ is compact,  its inverse is unbounded,
so a non-zero lower bound for  $\| A_{\bm_0} h_0  \|_{L^2(V)}$ has to be
considered rather than a lower bound for $\| h_0 \|_{H^1_0(R)}$.

\begin{thm}\label{main theorem}
Let $B$ be a closed and bounded subset of geometry parameters $\bm$ in $\RR^3$ such that 
  the distance condition \eqref{BProperty} holds. Let $A_\bm$ be the crack to boundary 
	operator defined by \eqref{AmDef}.
	Let $M_1$ and $M_2$ be two positive constants.
	There is a positive constant $C$ such that for all
   $\bm, \bm'$ in $B$, all $g$ in $H^1_0(R)$,
	and all $h$ in the set
	\bean \label{setS}
	S = \{ \varphi \in H^1_0(R): 
	\| \varphi \|_{H^1_0(R)} \leq M_2 \mbox{ and }
	M_1 \leq \| A_{\bm} \varphi \|_{L^2(V)} 
	\mbox{ for all } \bm \in B
	\},
\eean
	 the following estimate holds,
	\bean \label{main_est}
	\| A_\bm g -  A_{\bm'} h  \|_{L^2(V)}  \geq
 C  |\bm - \bm'|.
\eean
\end{thm}
\textbf{Proof:} 
Arguing by contradiction assume that there are sequences 
$\bm_n$ and $\bm_n'$ in $B$, $g_n$ in 
$H^1_0(R)$ and $h_n$ in $S$ such that $\bm_n \neq\bm_n '$ and
\bean
\label{to cont}
\lim_{n \to \infty} 
\f{\| A_{\bm_n} g_n -  A_{\bm_n'} h_n  \|_{L^2(V)}}{ |\bm_n - \bm_n '|}
= 0.
\eean
By compactness, we may assume that
$\bm_n$ converges to some $\tilde{\bm}$ in $B$,
$\bm_n' $ converges to some $\tilde{\tilde{\bm}}$ in $B$,
and $h_n$ converges weakly to some $\tilde{h}$ in $H^1_0(R)$.
Note that $A_{\bm_n}$ converges to $A_{\tilde{\bm}} $ in operator norm,
$A_{\bm_n'}$ converges to $A_{\tilde{\tilde{\bm}}} $ in operator norm,
and $A_{\bm_n'} h_n$ converges strongly to 
$A_{\tilde{\tilde{\bm}}}  \tilde{h}$ in $L^2(V)$. 
As $h_n$ is in $S$, $M_1 \leq \| A_{\bm_n'} h_n \|_{L^2(V)}$,
so at the limit we have that $\tilde{h} \neq 0$.
\\
Assume that $\tilde{\bm } \neq \tilde{\tilde{\bm}}$.
As,
\bea
\|( I - P_{\bm_n}) A_{\bm_n'} h_n  \|_{L^2(V)} \leq
\| A_{\bm_n} g_n -  A_{\bm_n'} h_n  \|_{L^2(V)},
\eea
where as earlier  $P_{\bm_n}$ is the orthogonal projection onto $\ov{{\cal R}(A_{\bm_n})}$,
$( I - P_{\bm_n}) A_{\bm_n'} h_n $ is also convergent to zero and at the limit 
we find $A_{\tilde{\tilde{\bm}}}  \tilde{h} =   P_{\tilde{\bm}} A_{\tilde{\tilde{\bm}}}  
\tilde{h} $,
so there is an $h'$ in $H^1_0(R)$ such that
$A_{\tilde{\tilde{\bm}}}  \tilde{h} =  A_{\tilde{\bm}}  h'
$. By uniqueness Theorem \ref{InverseProblemResult}, we infer that 
$\tilde{h} = 0$: contradiction. \\
We now know that $\tilde{\bm } =\tilde{\tilde{\bm}}$.
We first observe that,
\bean \label{first obs}
\lim_{n \to \infty} 
\f{\| ( I - P_{\bm_n}) (A_{\bm_n} -A_{\bm_n'}) h_n \|_{L^2(V)}}{ |\bm_n - \bm_n '|}
= 0.
\eean
Since $\f{\bm_n - \bm_n'}{|\bm_n - \bm_n'|}$ is a sequence on the unit sphere 
of $\RR^3$, after
possibly extracting a subsequence we may assume that it converges to some
$\bq$ with $|\bq|=1$.
Then as, 
\bea
A_{\bm_n} = A_{\tilde{\bm}} + (\bm_n - \tilde{\bm} ) \nabla A_{\tilde{\bm}}
+ |\bm_n - \tilde{\bm} | u_n, \\
A_{\bm_n'} = A_{\tilde{\bm}} + (\bm_n' - \tilde{\bm} ) \nabla A_{\tilde{\bm}}
+ |\bm_n' - \tilde{\bm} | v_n, 
\eea
where $\lim_{n \to \infty}  u_n = \lim_{n \to \infty} v_n =0$ in operator norm,
we infer that 
$
 \f{A_{\bm_n} -A_{\bm_n'}}{ |\bm_n - \bm_n '|}
$
converges to $\p_{\bq} A_{\tilde{\bm}}   $ in operator norm.
Since $h_n$ is weakly convergent to $\tilde{h}$ and the operator
$\p_{\bq} A_{\tilde{\bm}}   $ is compact, we obtain from the limit 
\eqref{first obs} the relation
$ (I - P_{\tilde{\bm}}) \p_{\bq} A_{\tilde{\bm}}  \tilde{h} =0$.  
This is akin to relation
\eqref{preg0} in the proof of Theorem \ref{stability2}. The same argument as in that 
proof will  show that $\tilde{h} =0$: contradiction. $\square$

\section{Conclusion and perspectives for future work}
We have shown in this paper a uniform stability result
for a planar crack inverse problem in half space governed by the Laplace equation.
We proved that reconstructing 
the plane containing the crack is Lipschitz stable despite the fact that the forcing term for the underlying PDE is unknown. This stability result holds under the assumption that the forcing
term $g$ is bounded above in  $H^1_0$ norm
and the Dirichlet data on the top boundary is bounded below in  $L^2$ norm.\\
In future work, we are planning to investigate 
 computational aspects of solving this crack
 inverse problem. In the case of the full   linear elasticity 
problem, this was done in \\
\cite{volkov2019stochastic} and in \cite{volkov2020stochastic}, where special random walk techniques were
designed and implemented.
Our research has shown
\cite{volkov2020stochastic} that finding 
  an  objective and automatic way of dealing with regularization parameters
	is possible thanks to random walks.
Conducting
 a similar study for 
the case governed by the
Laplace equation 
would present the advantage that
the related Green function is orders of magnitude
simpler to compute. 
It would therefore be
 easier to  focus on the development of 
novel reconstruction  algorithms for the 
posterior  probability distribution functions
of geometry parameters defining the crack,
unhindered  
 by the  cost of computing the 
Green function. 
 \vskip 15pt
\Large{\bf{Acknowledgments}} \\
\normalsize
This work was supported by
 Simons Foundation Collaboration Grant [351025].


\appendix
\section{A counter example for the unique solvability of the crack inverse problem in case of Lipschitz cracks}\label{counter-example}
\begin{figure}[htbp]
    \centering
      \includegraphics[scale=.5]{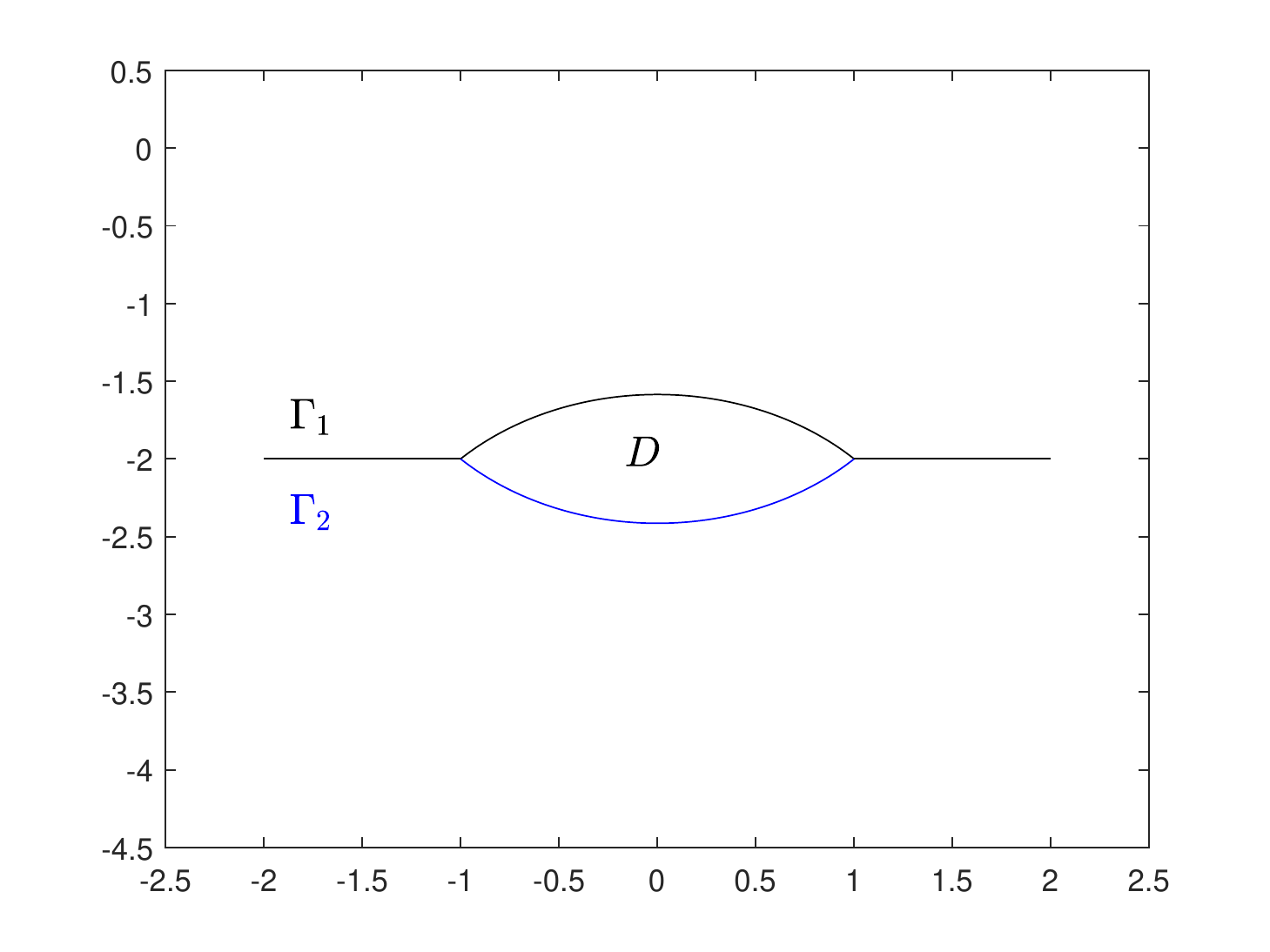}
         \caption{A cross section showing the two open surfaces $\Gamma_1$ and $\Gamma_2$ together
				with the domain $D$.}\label{counter}
\end{figure}

Define the open disk  in $\RR^2$,
\bea
R =\{  (x_1, x_2) \in \RR^2 : \s{x_1^2+ x_2^2} < 2 \}.
\eea
Define two Lipschitz functions from $R$ to $[-1 -\s{2},-3+ \s{2}]$,
\bea
\psi_1(x_1,x_2)= \left\{ 
\begin{array}{ll}
-2, &\mbox{ if } 1 \leq \s{x_1^2+ x_2^2} < 2\\
  -3 + \s{ 2 -x_1^2 -  x_2^2}, & \mbox{ if } 0 \leq \s{x_1^2+ x_2^2} <1,
	\end{array}
\right.
\eea
\bea
\psi_2(x_1,x_2)= \left\{ 
\begin{array}{ll}
-2, &\mbox{ if } 1 \leq \s{x_1^2+ x_2^2} < 2 \\
  -1 - \s{ 2 -x_1^2 -  x_2^2}, & \mbox{ if } 0 \leq \s{x_1^2+ x_2^2} <1.
	\end{array}
\right.
\eea
Cross sections of the graphs of $\psi_1$ and $\psi_2$ are sketched in 
Figure \ref{counter}.
Next, define the  function $g$ in $R$
\bea
g(x_1,x_2) = \left\{ 
\begin{array}{ll}
4 - x_1^2 - x_2^2, &\mbox{ if } 1 \leq  \s{x_1^2+ x_2^2} <  2\\
  3 , & \mbox{ if } 0 \leq \s{x_1^2+ x_2^2} <1.
	\end{array}
\right.
\eea
$g$ can be extended by zero outside $R$ to obtain a function in $H^1_0(\RR^2)$.
For $i=1,2$, let $\Gamma_i$ be the open surface defined by the graph of $\psi_i$  and define
$u^i$ be the function on $\RR^{3-}\sm \ov{\Gamma_i} $
\bea
u^i(\bx)=
\int_R H(x_1,x_2,x_3,y_1,y_2,\psi_i (y_1,y_2),\bn_i)
				g(y_1,y_2)\sigma_i
				dy_1dy_2, 
\eea
where $H$ is as in \eqref{Hdef},
$\sigma_i $ is defined almost everywhere on $R$ and equals 
$\s{1 + (\p_{x_1} \psi_i)^2+ (\p_{x_2} \psi_i)^2 } $,
and $\bn_i = ( -\p_{x_1} \psi_i, - \p_{x_2} \psi_i,1)/\sigma_i$, almost everywhere.
Let $u=u^1-u^2$ and $D$, sketched in Figure \ref{counter},  be the Lipshitz domain
defined by 
\bea
 -1 - \s{ 2 -x_1^2 -  x_2^2}   <x_3< -3 + \s{ 2 -x_1^2 -  x_2^2}, \quad
 \s{x_1^2+ x_2^2} <1.
\eea
Let $u=u^1-u^2$.
We observe that the exterior normal unit vector $\bn$ to $\p D$ is $\bn^1$ on $\p D \cap \Gamma_1$ and
 $-\bn^2$ on $\p D \cap \Gamma_2$, so  it follows that, 
\bea
u(\bx) = 
3 \int_{\p D} H(\bx, \by,\bn)
				d\sigma (\by),
\eea
for $\bx$ in $\RR^{3-} \sm \p D$.
By Green's theorem,  $u(\bx) = 0$ if $\bx \in  \RR^{3-} \sm \ov{D}$, 
$u(\bx) = -3 $ if $\bx \in D$.
We conclude that for all $(x_1, x_2)$ in $\RR^2$, $u^1(x_1, x_2,0) = u^2(x_1, x_2,0)$
and $\p_{x_3}u^1(x_1, x_2,0) = \p_{x_3}u^2(x_1, x_2,0)$, thus the crack inverse problem is not uniquely solvable
without additional assumptions.

\section{The geometry of $U_j$ for $j \geq 2$}\label{geom}
As $\psi_i: R_i \ri \RR$ is Lipschitz continuous, we can extend $\psi_i$ as a
Lipschitz continuous function from $\ov{R_i}$ to $ \RR$, $i=1,2$. 
\begin{lem} \label{intUj}
Let $j \geq 2$ and $(x_1,x_2,x_3)$ be in $U_j$. Then $(x_1,x_2) \in R_1 \cap R_2$
and $\psi_1(x_1,x_2)  \neq \psi_2 (x_1,x_2) $.
\end{lem}
\textbf{Proof:}
Arguing by contradiction assume that  $j \geq 2$, $(x_1,x_2,x_3) \in U_j$,
and $(x_1, x_2) \notin R_1$.\\
If $(x_1, x_2) \notin R_2$ then the half line
$\{(x_1,x_2, t x_3): t>0 \}$ is in $U$, so $(x_1,x_2,x_3) \in U_1$: contradiction. \\
It follows that $(x_1, x_2) \in R_2$. If $\psi_2(x_1,x_2) \geq x_3$ then the half-line
$\{(x_1,x_2, t x_3): t>1 \}$ is in $U_1$. As $U_j $ is open, for some $\alpha >1$,
$(x_1,x_2, \alpha x_3) \in U_j$: contradiction. 
If  $\psi_2(x_1,x_2) < x_3$, then the line segment $\{(x_1,x_2, t x_3): t \in (0,1] \}$
is in $U_1$: contradiction. 
At this stage we conclude that $(x_1, x_2) \in R_1$ and similarly $(x_1, x_2) \in R_2$. There remains to show that $\psi_1(x_1,x_2)  \neq \psi_2 (x_1,x_2) $.
Arguing by contradiction, assume that $\psi_1(x_1,x_2)  =\psi_2 (x_1,x_2) $.
Then the half-line $L_1 = \{(x_1,x_2, t \psi_1(x_1,x_2)): t>1 \}$
and the line segment $L_2 = \{(x_1,x_2, t \psi_1(x_1,x_2)): t \in (0,1) \} $ are in $U_1$.
Since  $j \geq 2$ and $(x_1,x_2,x_3) \in U_j$, $(x_1, x_2, x_3)$ 
is neither in $L_1$ nor in $L_2$. It follows that $x_3= \psi_1(x_1,x_2) $.
But as $U_j$ is open $(x_1, x_2, \alpha x_3) \in U_j$ for some $\alpha >0$: contradiction. 
$\square$
\begin{lem} \label{pos or neg}
Fix $j \geq 2$.
Then  for all  $(x_1,x_2,x_3)$  in $U_j$, 
 $\psi_1(x_1,x_2) < \psi_2 (x_1,x_2) $, or 
 for all  $(x_1,x_2,x_3)$  in $U_j$, 
 $\psi_1(x_1,x_2) > \psi_2 (x_1,x_2) $ 
\end{lem}
\textbf{Proof:}
The function $\phi_j: U_j \ri \RR$, $\phi_j(x_1,x_2,x_3) = \psi_1(x_1,x_2) - \psi_2 (x_1,x_2)$,
is continuous. As $U_j$ is connected, $\phi_j (U_j)$ is connected in $\RR$.
By Lemma  \ref{intUj}, $0 \notin \phi_j (U_j)$.
$\square$
\begin{lem} \label{omega j}
Let $j \geq 2$ and 
$\Omega_j = \{ (x_1, x_2) \in \RR^2: (x_1,x_2,x_3) \in U_j\}$.
Then $\Omega_j$ is a  bounded and connected open subset of $\RR^2$ 
and for all $ (x_1, x_2) $ in $\p \Omega_j $, $\psi_1(x_1,x_2) = \psi_2(x_1,x_2)$.
\end{lem}
\textbf{Proof:}
The first assertion is clear since $U_j$ is open, bounded, and connected in $\RR^3$.
Arguing by contradiction, assume that for some 
$ (x_1, x_2) $ in $\p \Omega_j $, $\psi_1(x_1,x_2) \neq \psi_2(x_1,x_2)$.
Without loss of generality we may assume that $\psi_1(x_1,x_2) < \psi_2(x_1,x_2)$.
Since $R_1 \cap R_2$ is open, $ (x_1, x_2) $ is in $R_1 \cap R_2$,
or it is on its boundary, $\p (R_1 \cap R_2)$.
In the first case, by continuity, there is a positive $\epsilon$ and a positive $\alpha$
such that if $(y_1, y_2)$ is in $(x_1 - \alpha, x_1 + \alpha) \times
 (x_2 - \alpha, x_2 + \alpha) \subset (R_1 \cap R_2)$, 
$\psi_1(y_1,y_2) < \psi_2(y_1,y_2) + \epsilon$.
But then for any element $(z_1, z_2,z_3)$ in $(x_1 - \alpha, x_1 + \alpha) \times
  (x_2 - \alpha, x_2 + \alpha) \times 
	(\psi_1(x_1,x_2) , \psi_2(x_1,x_2))$,  the line segment from $(z_1,z_2, z_3)$ to
	$(x_1, x_2, \psi_1(x_1,x_2) + \f{\epsilon}{2})$
is included in $U$, thus the whole set
$(x_1 - \alpha, x_1 + \alpha) \times
 (x_2 - \alpha, x_2 + \alpha) \times (\psi_1(x_1,x_2) , \psi_2(x_1,x_2))$
is in the same connected component $U_j$: this contradicts that $ (x_1, x_2) $ is
 in $\p \Omega_j $.\\
In the second case, $ (x_1, x_2) $ is in $\p (R_1 \cap R_2)$. Without loss of generality say that
it is in $\p R_1 $. By a continuity argument, a small open ball centered at 
$(x_1, x_2, \psi_1(x_1,x_2) + \f{\epsilon}{2})$ for some small $\epsilon >0$ is included in $U_1$
and intersects $U_j$, which is a contradiction. $\square$ \\\\
Using lemmas  \ref{intUj}, \ref{pos or neg},  and \ref{omega j}, we see that $U_j$ is simple in the 
$x_3$ direction and that for any integrable function $f$ in $U_j$,
\bean \label{int form Uj}
\int_{U_j} f = \int_{\Omega_j} \int_{\psi_k(x_1,x_2)}^{\psi_l(x_1,x_2)}
f(x_1,x_2,x_3) dx_3 dx_1 d x_2,
\eean
where depending on $j \geq 2$, $(k,l) = (1,2)$ or $(k,l) = (2,1)$.

\section{Proof of Lemma \ref{thirdjumplem}}\label{appen prove jump}
We first state and prove two intermediary results.
\begin{lem} \label{firstjumplem}
    Let $\Gamma$ be an open surface in $\doubleR^3$ included in the plane $x_3=0$ and let $g$ be in $C_c^\infty(\Gamma)$. 
		Using the notation for $\bx=(x_1,x_2,0)$ in $\Gamma$,
		\bea
		 [ v(\bx)] = \lim_{\epsilon \ri 0^+}
		v(x_1,x_2, \epsilon) - v(x_1,x_2, -\epsilon),
		\eea
		 the following jump formulas across $\Gamma$ hold for $i=1, 2$,
\begin{equation}
    \label{jumpformula1}
    [\int_\Gamma G(\bx,\by,\bev_3)g(\by)d\by]=g(\bx),
\end{equation}
\begin{equation}
    \label{jumpformula2}
    [\int_\Gamma (\partial_{y_i}G)(\bx,\by,\bev_3)g(\by)d\by]=-\partial_{x_i}g(\bx),
\end{equation}
\begin{equation}
    \label{jumpformula3}
    [\int_\Gamma G(\bx,\by,\bev_i)g(\by)d\by]=0,
\end{equation}
\begin{equation}
    \label{jumpformula4}
    [\int_\Gamma (\partial_{y_3}G)(\bx,\by,\bev_3)g(\by)d\by]=0,
\end{equation}

\begin{equation}
    \label{jumpformula42}
    [\p_{x_3} \int_\Gamma G(\bx,\by,\bev_3)g(\by)d\by]=0,
\end{equation}

\begin{equation}
    \label{jumpformula5}
    [\p_{x_3}\int_\Gamma G(\bx,\by,\bev_i)g(\by)d\by]= \partial_{x_i}g(\bx),
\end{equation}

\begin{equation}
    \label{jumpformula6}
    [\p_{x_3}\int_\Gamma (\partial_{y_3}G)(\bx,\by,\bev_3)g(\by)d\by]= \Delta g(\bx).
\end{equation}

\end{lem}
\textbf{Proof:}
We can find a domain $D$ of class $C^2$ such
that $\Gamma \subset \p D$ and the exterior normal vector to 
$\p D$ is equal to $\bev_3$ on $\Gamma$. Extend  $g$  by zero on 
$\p D\setminus \Gamma$ to obtain a function in $C^2(\p D)$.
Now, given the definition of $G$ \eqref{Gdef},
\eqref{jumpformula1} is just the classical jump formula for the double layer potential.
To show  \eqref{jumpformula2}, we note that for $x_3 \neq 0$, $\by$ in $\Gamma$,
 and $i=1,2$,
\bean \int_\Gamma (\partial_{y_i}G)(\bx,\by,\bev_3)g(\by)d\by
= - \int_\Gamma G(\bx,\by,\bev_3) \p_i g(\by)d\by,
\eean
and we just apply \eqref{jumpformula1}.
Similarly, to show \eqref{jumpformula3}, we note that for $x_3 \neq 0$,
$\by$ in $\Gamma$, and $i=1,2$,
\bea \int_\Gamma G(\bx,\by,\bev_i)g(\by)d\by
= - \int_\Gamma \Phi(\bx,\by) \p_i g(\by)d\by,
\eea
and we use the well known continuity of the single layer potential.\\
 Identity \eqref{jumpformula42} is due to the well known continuity property of the normal derivative
of the double layer potential. \\
To show \eqref{jumpformula4} and \eqref{jumpformula5}, we note that for $x_3 \neq 0$ and $\by$ in $\Gamma$,
since $\Phi$ depends on $\bx$ and $\by$ only through $\bx - \by$, for $i=1,2$ or $3$,
\bea
\p_{x_3} G(\bx,\by,\bev_i) = \p_{x_3} (\nabla_\by \Phi (\bx,\by) \cdot \bev_i )
= -  \nabla_\by  (\p_{y_3}\Phi (\bx,\by)) \cdot \bev_i .
\eea
If $i=3$, \eqref{jumpformula4} is now clear due to \eqref{jumpformula42}. 
To show \eqref{jumpformula5}, 
for $i=1,2$,
\bea
\p_{x_3}\int_\Gamma G(\bx,\by,\bev_i)g(\by)d\by 
= \int_\Gamma (\p_{y_3}\Phi (\bx,\by))\partial_{y_i}g(\by)   d\by  
= \int_\Gamma   G(\bx,\by,\bev_3) \partial_{y_i}g(\by)   d\by,  
\eea
and we can apply  \eqref{jumpformula1}.
To show \eqref{jumpformula6}, 
we use again that $\Phi$ depends on $\bx$ and $\by$ only through $\bx - \by$,
so if  $x_3 \neq 0$ and $\by$ is  in $\Gamma$,
\bea
\p_{x_3} \p_{y_3} G(\bx,\by,\bev_3) = -  \p_{y_3}^2 G(\bx,\by,\bev_3) =
(\p_{y_1}^2 + \p_{y_2}^2) G(\bx,\by,\bev_3) ,
\eea
where we have used that $\Delta_y \Phi (\bx, \by) =0$ if $\bx \neq \by$.
Now as $g$ is in $C^\infty_c(\Gamma)$,
\bea
\int_\Gamma (\p_{y_1}^2 + \p_{y_2}^2) G(\bx,\by,\bev_3)g(\by)d\by 
= \int_\Gamma G(\bx,\by,\bev_3) (\p_{y_1}^2 + \p_{y_2}^2) g(\by)   d\by ,
\eea
and the result follows from \eqref{jumpformula1}.

\begin{lem} \label{secondjumplem}
    Let $\Gamma$ be as in lemma \ref{firstjumplem} and $g$ be in 
		$L^2(\Gamma)$. 
		Then the jump
		formulas (\ref{jumpformula1}-\ref{jumpformula6}) 
		hold in a weak sense, that is for any $\phi$ in $C^\infty_c(\Gamma)$ 
		and $i=1, 2$,
		\begin{equation}
    \label{jumpformula1w}
    \int_\Gamma [\int_\Gamma G(\bx,\by,\bev_3)g(\by) d  \by ] \phi (\bx) d\bx =
		\int_\Gamma g(\bx)  \phi (\bx) d  \bx,
\end{equation}
\begin{equation}
    \label{jumpformula2w}
    \int_\Gamma[  \int_\Gamma (\partial_{y_i}G)(\bx,\by,\bev_3)g(\by) d  \by] \phi (\bx) d\bx
		=\int_\Gamma g(\bx)  \partial_{x_i}\phi (\bx) d  \bx,
\end{equation}
\begin{equation}
    \label{jumpformula3w}
    \int_\Gamma [\int_\Gamma G(\bx,\by,\bev_i)g(\by)d  \by] \phi (\bx)  d\bx =0,
\end{equation}
\begin{equation}
    \label{jumpformula4w}
    \int_\Gamma [\int_\Gamma (\partial_{y_3}G)(\bx,\by,\bev_3)d  \by] \phi (\bx)  d\bx =0,
\end{equation}

\begin{equation}
    \label{jumpformula42w}
    \int_\Gamma [\p_{x_3} \int_\Gamma G(\bx,\by,\bev_3)d\by ] \phi (\bx) d  \bx=0,
\end{equation}

\begin{equation}
    \label{jumpformula5w}
    \int_\Gamma[\p_{x_3}\int_\Gamma G(\bx,\by,\bev_i)d\by ] \phi (\bx) d  \bx  =
		- \int_\Gamma g(\bx)  \partial_{x_i} \phi (\bx) d \bx ,
\end{equation}

\begin{equation}
    \label{jumpformula6w}
    \int_\Gamma[\p_{x_3}\int_\Gamma (\partial_{y_3}G)(\bx,\by,\bev_3)
		g(\by)d\by]\phi (\bx) d  \bx
		= \int_\Gamma g(\bx)  \Delta \phi (\bx) d \bx.
\end{equation}

\end{lem}
\textbf{Proof:}
Let $\epsilon$ be positive and set,
\bea
\bx^\pm = (x_1, x_2, \pm \epsilon), \,  \by^\pm = (y_1, y_2, \pm \epsilon).
\eea
As $\Phi$ depends on $\bx$ and $\by$ only through $\bx - \by$,
\bea
&G(\bx^+,\by, \bev_3) - G(\bx^-,\by,\bev_3) \\
= & \nabla_\by\Phi(\bx^+, \by) \cdot \bev_3  - \nabla_\by\Phi(\bx^-, \by) \cdot \bev_3 \\
= & -\nabla_\bx\Phi(\bx^+, \by) \cdot \bev_3  + \nabla_\bx\Phi(\bx^-, \by) \cdot \bev_3 \\
= & -\nabla_\bx\Phi(\bx, \by^-) \cdot \bev_3  + \nabla_\bx\Phi(\bx, \by^+) \cdot \bev_3,
\eea
thus by \eqref{jumpformula1},
\bean \label{uni}
\lim_{\epsilon \ri 0 } \int_\Gamma (G(\bx^+,\by, \bev_3) - G(\bx^-,\by,\bev_3)) \phi(\bx)
 d \bx = \phi (\by).
\eean
But the convergence in \eqref{uni} is known to be uniform for $\by$ in $\Gamma$,
 thus \eqref{jumpformula1w} follows.\\
To prove \eqref{jumpformula2w}, we observe that
 \bea
\int_\Gamma (\partial_{y_i}G)(\bx^\pm,\by,\bev_3)\phi (\bx) d  \bx  
= - \int_\Gamma (\partial_{x_i}G)(\bx^\pm,\by,\bev_3)\phi (\bx) d  \bx 
= \int_\Gamma G(\bx^\pm,\by,\bev_3) \partial_{x_i}\phi (\bx) d  \bx 
\eea
so the result follows from \eqref{jumpformula1w}. \\
Identities (\ref{jumpformula3w}-\ref{jumpformula5w}) can be derived likewise from lemma \ref{firstjumplem}.
We now show the derivation of \eqref{jumpformula6w}.
As,
\bea
\int_\Gamma \p_{x_3} \p_{y_3} G(\bx^\pm, \by, \bev_3 ) \phi(\bx) d\bx 
= \int_\Gamma -\p_{x_3}^2  G(\bx^\pm, \by, \bev_3 ) \phi(\bx) d\bx \\
= \int_\Gamma (\p_{x_1}^2 +   \p_{x_2}^2 )G(\bx^\pm, \by, \bev_3 ) \phi(\bx) d\bx 
= \int_\Gamma G(\bx^\pm, \by, \bev_3 ) \Delta \phi(\bx) d\bx, 
\eea
the result follows from \eqref{jumpformula1w}. 
$\square$\\

To prove Lemma \ref{thirdjumplem}, it now suffices to
 extend Lemma \ref{secondjumplem} to the case of planar open surfaces
that are not necessarily included in the plane $x_3=0$.
Since the fundamental solution $\Phi$ given by \eqref{fund}
satisfies
$
\Phi(T\bx,T\by)=\Phi(\bx,\by),
$
for all $\bx$ and $\by$ in $\doubleR^3$ such that $\bx \neq \by$ 
where $T$   is any 
rotation   or translation
of $\doubleR^3$, we can obtain formulas
\eqref{jumpformula31} to \eqref{jumpformula36} by
 a straightforward generalization of Lemma
\ref{secondjumplem}.
   Formula \eqref{jumpformula343}
 can be obtained  from \eqref{jumpformula342}
by applying Green's formula on $\Gamma$.
$\square$

\bibliography{ref}{}

\begin{thebibliography}{10}

\bibitem{aspri2020arxiv}
Andrea Aspri, Elena Beretta, and Anna~L Mazzucato.
\newblock Dislocations in a layered elastic medium with applications to fault
  detection.
\newblock {\em preprint arXiv:2004.00321v1}, 2020.

\bibitem{aspri2020analysis}
Andrea Aspri, Elena Beretta, Anna~L Mazzucato, and V~Maarten.
\newblock Analysis of a model of elastic dislocations in geophysics.
\newblock {\em Archive for Rational Mechanics and Analysis}, 236(1):71--111,
  2020.

\bibitem{costabel1988boundary}
Martin Costabel.
\newblock Boundary integral operators on lipschitz domains: elementary results.
\newblock {\em SIAM Journal on Mathematical Analysis}, 19(3):613--626, 1988.

\bibitem{dascalu2000fault}
Cristian Dascalu, Ioan~R Ionescu, and Michel Campillo.
\newblock Fault finiteness and initiation of dynamic shear instability.
\newblock {\em Earth and Planetary Science Letters}, 177(3):163--176, 2000.

\bibitem{ding1996proof}
Zhonghai Ding.
\newblock A proof of the trace theorem of sobolev spaces on lipschitz domains.
\newblock {\em Proceedings of the American Mathematical Society},
  124(2):591--600, 1996.

\bibitem{friedman1989determining}
Avner Friedman and Michael Vogelius.
\newblock Determining cracks by boundary measurements, 1989.
\newblock http://conservancy.umn.edu/bitstream/handle/11299/4926/476.pdf.

\bibitem{gagliardo1957caratterizzazioni}
Emilio Gagliardo.
\newblock Caratterizzazioni delle tracce sulla frontiera relative ad alcune
  classi di funzioni in $ n $ variabili.
\newblock {\em Rendiconti del seminario matematico della universita di Padova},
  27:284--305, 1957.

\bibitem{ionescu2006inverse}
Ioan~R Ionescu and Darko Volkov.
\newblock An inverse problem for the recovery of active faults from surface
  observations.
\newblock {\em Inverse problems}, 22(6):2103, 2006.

\bibitem{ionescu2008earth}
Ioan~R Ionescu and Darko Volkov.
\newblock Earth surface effects on active faults: An eigenvalue asymptotic
  analysis.
\newblock {\em Journal of Computational and Applied Mathematics},
  220(1):143--162, 2008.

\bibitem{kato2013perturbation}
Tosio Kato.
\newblock {\em Perturbation theory for linear operators}, volume 132.
\newblock Springer Science \& Business Media, 2013.

\bibitem{triki2019stability}
Faouzi Triki and Darko Volkov.
\newblock Stability estimates for the fault inverse problem.
\newblock {\em Inverse problems}, 35(7), 2019.

\bibitem{volkov2017determining}
D.~Volkov, C.~Voisin, and Ionescu I.R.
\newblock Determining fault geometries from surface displacements.
\newblock {\em Pure and Applied Geophysics}, 174(4):1659--1678, 2017.

\bibitem{volkov2020stochastic}
Darko Volkov.
\newblock A stochastic approach to mixed linear and nonlinear inverse problems
  with applications to seismology.
\newblock {\em preprint arXiv:2007.05347}, 2020.

\bibitem{volkov2019stochastic}
Darko Volkov and Joan~Calafell Sandiumenge.
\newblock A stochastic approach to reconstruction of faults in elastic half
  space.
\newblock {\em Inverse Problems \& Imaging}, 13(3):479--511, 2019.

\bibitem{volkov2017reconstruction}
Darko Volkov, Christophe Voisin, and Ioan Ionescu.
\newblock Reconstruction of faults in elastic half space from surface
  measurements.
\newblock {\em Inverse Problems}, 33(5), 2017.

\end{thebibliography}
\bibliographystyle{plain}

\end{document}